\documentclass[11pt,british,letterpaper]{amsart}

\newif\ifpdf
\ifx\pdfoutput\undefined\pdffalse\else\pdfoutput=1\pdftrue\fi\newif\pdf
\ifpdf\relax\else
\usepackage[T1]{fontenc}
\fi

\usepackage{babel,eucal,url,amssymb,stmaryrd,
enumerate,amscd,
}

\usepackage[pagebackref]{hyperref}

\theoremstyle{plain}
\newtheorem{thm}{Theorem}[section]

\newtheorem{pro}[thm]{Proposition}
\newtheorem{co}[thm]{Corollary}

\theoremstyle{definition}

\theoremstyle{remark}

\newcommand{\Gtwo}{\ifmmode{{\rm G}_2}\else{${\rm G}_2$}\fi}

\newcommand{\LC}{{\nabla^g}}

\def\sideremark#1{\ifvmode\leavevmode\fi\vadjust{\vbox to0pt{\vss
 \hbox to 0pt{\hskip\hsize\hskip1em
 \vbox{\hsize2.5cm\tiny\raggedright\pretolerance10000
 \noindent #1\hfill}\hss}\vbox to8pt{\vfil}\vss}}}%

\date{\today}
\begin{document}
\title[Twistors of Almost Quaternionic
Manifolds]%
{Twistors  of Almost Quaternionic Manifolds}
\date{\today}

\author{Stefan Ivanov}
\address[Ivanov, Minchev, Zamkovoy]{University of Sofia "St. Kl. Ohridski"\\
Faculty of Mathematics and Informatics,\\ Blvd. James Bourchier
5,\\ 1164 Sofia, Bulgaria} \email{ivanovsp@fmi.uni-sofia.bg,
minchev@fmi.uni-sofia.bg, zamkovoy@fmi.uni-sofia.bg}

\author{Ivan Minchev}

\author{Simeon Zamkovoy}

\thanks{Partially supported by a Contract 154/2005
with the University of Sofia "St. Kl. Ohridski"}

\begin{abstract}
We investigate the integrability of  almost complex structures on
the twistor space of an almost quaternionic manifold constructed
with the help of a quaternionic connection. We show that if there
is an integrable structure it is independent on the quaternionic
connection. In dimension four, we express the anti-self-duality
condition in terms of the Riemannian Ricci forms with respect to
the associated quaternionic structure.
\\[8mm] {\bf Keywords.}  Almost Quaternionic, Torsion, Twistors, QKT-spaces.
\\[8mm]
{\bf AMS Subject Classification:}  Primary 53C25, Secondary
53C15, 53C56, 32L25, 57S25
\end{abstract}

\maketitle

\setcounter{tocdepth}{2} \tableofcontents

\section{Introduction and statement of the results}

An almost hyper-complex structure on a 4n-dimensional manifold $M$
is a triple $H=(J_{\alpha}), \alpha=1,2,3$, of almost complex
structures $J_{\alpha}:TM\rightarrow TM$ satisfying the
quaternionic identities $J_{\alpha}^2=-id$ and
$J_1J_2=-J_2J_1=J_3$. When each $J_{\alpha}$ is a complex
structure, $H$ is said to be a hyper complex structure on $M$.

An almost quaternionic structure on $M$ is a rank-3 subbundle
$\mathcal Q \subset End(TM)$ which is locally spanned by an almost
hyper-complex structure $H=(J_{\alpha})$. Such a locally defined
triple $H$ is called an admissible basis of $\mathcal Q$. A linear
connection $\nabla$ on $TM$ is called a quaternionic connection if
$\nabla$ preserves $\mathcal Q$, i.e. $\nabla_X\sigma\in \Gamma
(\mathcal Q)$ for all vector fields $X$ and smooth sections
$\sigma \in \Gamma (\mathcal Q)$. An almost quaternionic structure
is said to be quaternionic if there is a torsion-free quaternionic
connection. A $\mathcal Q$-hermitian metric is a Riemannian metric
which is Hermitian with respect to each almost complex structure
in $\mathcal Q$. An almost quaternionic (resp. quaternionic)
manifold with $\mathcal Q$-hermitian metric is called an almost
quaternionic Hermitian (resp. quaternionic hermitian) manifold

For $n=1$ an almost quaternionic structure is the same as an
oriented conformal structure and it turns out to be always
quaternionic.  The existence of a (local) hyper-complex structure
is a strong condition since the integrability of the (local)
almost hyper-complex structure implies that the corresponding
conformal structure is anti-self-dual \cite{2}.

When $n\ge 2$, the existence of torsion-free quaternionic
connection is a strong condition which is equivalent to the
1-integrability of the associated $GL(n,{\bf H})Sp(1)$ structure
\cite{Bon,Ob,Sal4}.  If the Levi-Civita connection $\LC$ of a
quaternionic hermitian manifold $(M,g,\mathcal Q)$ is a
quaternionic connection then $(M,g,\mathcal Q)$ is called a
Quaternionic K\"ahler manifold (briefly QK manifold). This
condition is equivalent to the statement that the holonomy group
of $g$ is contained in $Sp(n)Sp(1)$ \cite{A1,A2,S1,S2,Ish}.  If on
a QK manifold there exists an admissible basis $(H)$ such that
each almost complex structure $(J_{\alpha})\in(H), \alpha =1,2,3$
is parallel with respect to the Levi-Civita connection then the
manifold is called hyperK\"ahler (briefly HK). In this case the
holonomy group of $g$ is contained in Sp(n).

The various notions of quaternionic manifolds arise in a natural
way from the theory of supersymmetric sigma models as well as in
string theory.  The geometry of the target space of
two-dimensional sigma models with extended supersymmetry is
described by the properties of a metric connection with torsion
\cite{HP1,HP2}.  The geometry of (4,0) supersymmetric
two-dimensional sigma models without Wess-Zumino term (torsion) is
a hyperK\"ahler manifold. In the presence of torsion the geometry
of the target space becomes hyperK\"ahler with torsion (briefly
HKT) \cite{HP3}. This means that the complex structures
$J_{\alpha}, \alpha=1,2,3$, are parallel with respect to a metric
quaternionic connection with totally skew-symmetric torsion
\cite{HP3}. Local (4,0) supersymmetry requires that the target
space of two dimensional sigma models with Wess-Zumino term be
either HKT or quaternionic K\"ahler with torsion (briefly QKT)
\cite{Hish} which means that the quaternionic subbundle is
parallel with respect to a metric linear connection with totally
skew-symmetric torsion and the torsion 3-form is of type
(1,2)+(2,1) with respect to all almost complex structures in
$\mathcal Q$. The target space of two-dimensional (4,0)
supersymmetric sigma models with torsion coupled to (4,0)
supergravity is a QKT manifold \cite{HOP}.

The main object of interest in this article are properties of the
twistor space of an almost quaternionic manifold. The twistor
space $Z$ is the unit sphere bundle with fibre $S^2$ consisting of
all almost complex structures compatible with the given almost
quaternionic structure $\mathcal Q$. We consider two almost
complex structures $I^{\nabla}_1$, $I^{\nabla}_2$ on the twistor
space $Z$ over an almost quaternionic manifold generated by a
quaternionic connection $\nabla$. The structure $I^{\LC}_{1}$ was
originally constructed for four-dimensional Riemannian manifold in
\cite{2}, $I^{\LC}_{2}$ is described in \cite{11}. In both cases
the horizontal space of the Levi-Civita connection is used. It is
shown in \cite{2} that $I^{\LC}_1$ is integrable exactly when the
conformal structure is ant-self-dual while in \cite{11} the
non-integrability of $I^{\LC}_2$ is proved. The integrability of
$I^{\LC}_1$ for quternionic K\"ahler space is established in
\cite{S1} and then it is generalized for quaternionic manifold in
\cite{S3}. For HKT and QKT spaces the integrability of
$I^{\nabla}_1$ is proved in \cite{HP3} and \cite{HOP},
respectively. Hermitian and K\"ahler geometry of the twistor space
is investigated in \cite{H,S1}. The almost hermitian geometry of
the twistor space over a 4-dimensional Riemannian manifold is
studied in \cite{DM}. The almost hermitian geometry of the twistor
space over a quaternionic K\"ahler manifold, HKT and QKT spaces is
considered in \cite{AGI,IM}.

In the present note we investigate the dependence on the
quaternionic connection of these naturally defined almost complex
structures on the twistor space over an almost quaternionic
manifold. We obtain conditions on the quaternionic connection
which imply the coincidence of the corresponding almost complex
structures (Corollary~\ref{t2.7}, Corollary~\ref{t2.72}). We show
that the existence of an integrable almost complex structure on
the twistor space does not depend on the quaternionic connection
and it is equivalent to the condition that the almost quaternionic
$4n$-manifold is quaternionic for $n\ge2$ (Theorem~\ref{t2.1},
Theorem~\ref{main}). In particular, we show that the integrable
almost complex structures $I_1^{\LC}$ and $I_1^{\nabla}$ on a QKT
space coincide.

In dimension four we find new relations between the Riemannian
Ricci forms, i.e. the 2-forms which determine the
$Sp(1)$-component of the Riemannian curvature, which are
equivalent to the anti-self-duality of the oriented conformal
structure corresponding to a given quaternionic structure
(Theorem~\ref{four}).

\section{Preliminaries}

Let ${\bf H}$ be the quaternions and identify ${\bf H}^{n}={\bf
R}^{4n}$. To fix notation we assume that ${\bf H}$ acts on ${\bf
H}^{n}$ by right multiplication. This defines an antihomomorphism
\begin{gather*}\lambda :\{{\rm
unit \:quaternions}\} =\\= \{x+j_1y+j_2z+j_3w\ |\
x^2+y^2+z^2+w^2=1\} \longrightarrow SO(4n)\subset GL(4n,{\bf
R}),\end{gather*} where our convention is that $SO(4n)$ acts on
${\bf H}^{n}$ on the left. Denote the image by $Sp(1)$ and let
$J^0_{1}= -\lambda (j_1), J^0_{2}=-\lambda (j_2), J^0_{3}=-\lambda
(j_3)$. The Lie algebra of $Sp(1)$ is
$sp(1)=span\{J^0_{1},J^0_{2},J^0_{3}\}$ and we have
\begin{equation}\nonumber
{J^0_1}^2={J^0_2}^2={J^0_3}^2=-1,\qquad
J^0_1J^0_2=-J^0_2J^0_1=J^0_3.
\end{equation}
Define $GL(n,H)=\{A\in GL(4n,{\bf R}): A(sp(1))A^{-1}=sp(1) \}$.
The Lie algebra of $GL(n,H)$ is $gl(n,H)=\{A\in gl(4n,{\bf R}):
AB=BA$ for all $B\in sp(1)\}$.

Let $(M,\mathcal Q)$ be an almost quaternionic manifold and
$H=(J_a), a=1,2,3$ be an admissible local basis. Let $B\in
\Lambda^2(TM)$. We say that $B$ is of type $(0,2)_{J_a}$ with
respect to $J_a$ if $$B(J_aX,Y)=-J_aB(X,Y).$$ We denote this space
by $\Lambda_{J_a}^{0,2}$. The projection $B^{0,2}_{J_a}$ is given
by
\begin{gather*}
B^{0,2}_{J_a}(X,Y)=\frac{1}{4}
\left(B(J_aX,J_aY)-B(X,Y)-J_aB(J_aX,Y)-J_aB(X,J_aY)\right).
\end{gather*}
For example, the Nijenhuis tensor $N_a\in \Lambda^{0,2}_{J_a}$.

We denote the space of quaternionic connections on an almost
quaternionic manifold by $\Delta(\mathcal Q)$.

Let $\nabla\in\Delta(\mathcal Q)$ be a quaternionic connection on
an almost quaternionic manifold $(M,\mathcal Q)$. This means that
there exist locally defined 1-forms $\omega_{\alpha}, \alpha
=1,2,3$ such that
\begin{gather}\label{zzv}
 \nabla J_a=-\omega_b\otimes J_c+\omega_c\otimes J_b.
\end{gather}
Here and further $(a,b,c)$ stands for a cyclic permutation of
$(1,2,3)$.

It follows from \eqref{zzv} that the curvature $R^{\nabla}$ of any
quaternionic connection $\nabla\in\Delta(\mathcal Q)$ satisfies
the relations
\begin{gather}\label{rel1}
[R^{\nabla},J_a]=-A_b\otimes J_c+A_c\otimes J_b, \qquad
A_a=d\omega_a+\omega_b\wedge\omega_c.
\end{gather}
The Ricci 2-forms of a quaternionic connection are defined by
\begin{gather*}
 \rho^{\nabla}_a(X,Y)=-\frac12Tr\left(Z\longrightarrow
J_aR^{\nabla}(X,Y)Z\right), \quad a=1,2,3.
\end{gather*}
It is easy to see, using \eqref{rel1}, that the Ricci forms are
given by  $$\rho^{\nabla}_a=d\omega_a+\omega_b\wedge\omega_c. $$
We split the curvature of $\nabla$ into $gl(n,H)$-valued part
$(R^{\nabla})'$ and $sp(1)$-valued part $(R^{\nabla})''$ following
the classical scheme (see e.g. \cite{AM,Ish,Bes})
\begin{pro}\label{p1}
The curvature of  an almost quaternionic
connection on $M$ splits as follows
\begin{eqnarray}\nonumber
R^{\nabla}(X,Y)&=&(R^{\nabla})'(X,Y)+{1\over
  2n}(\rho^{\nabla}_1(X,Y)J_1+\rho^{\nabla}_2(X,Y)J_2+\rho^{\nabla}_3(X,Y)J_3),\\ \nonumber & &
  [(R^{\nabla})'(X,Y),J_a]=0,\ \ \ a=1,2,3.
\end{eqnarray}
\end{pro}
Let $\Omega, \Theta$ be the curvature 2-form and the torsion
2-form of $\nabla$ on the principal $GL(n,H)Sp(1)$-bundle
$\mathcal Q(M)$, respectively (\cite{15}). We denote the splitting
of the $gl(n,H)\oplus sp(1)$-valued curvature 2-form $\Omega$  on
$\mathcal Q(M)$ according to Proposition~\ref{p1}, by $\Omega
=\Omega' +\Omega''$, where $\Omega' $ is a $gl(n,H)$-valued
 2-form and  $\Omega''$  is  a
$sp(1)$-valued form.  Explicitly,  $$ \Omega ''= \Omega
''_{1}J^0_1 + \Omega''_{2}J^0_2 + \Omega''_{3}J^0_3,$$ where
$\Omega''_a, a=1,2,3$, are 2-forms.  If $\xi, \eta, \zeta \in {\bf
R^{4n}}$, then the 2-forms $\Omega''_a, a=1,2,3$, are given by
\begin{equation}\label{e1}
\Omega''_a(B(\xi ),B(\eta
))=\frac{1}{2n}\rho_a(X,Y), \quad X=u(\xi), Y=u(\eta).
\end{equation}

\section{Twistor  space of almost quaternionic manifolds}

In this section we adapt the setup from \cite{19,7} to incorporate
a torsion. Our discussion is very close to that of \cite{AGI,IM}.

Let $M$ be a $4n$-dimensional manifold endowed with  an almost
quaternionic structure $\mathcal Q$. Let $J_1,J_2, J_3$ be an
admissible basis of $\mathcal Q$ defined in some neighborhood of a
given point $p\in M$. Any linear frame $u$ of $T_pM$ can be
considered as an isomorphism $u:{\bf R}^{4n}\longrightarrow T_pM.$
If we pick such a frame $u$ we can define a subspace of the space
of the all endomorphisms of $T_pM$ by $u(sp(1))u^{-1}.$ Clearly,
this subset is a quaternionic structure at the point $p$ and in
the general case this quaternionic structure is different from
${\mathcal Q}_p.$ We define $\mathcal Q(M)$ to be the set of all
linear frames $u$ which satisfy $u(sp(1))u^{-1}={\mathcal Q}.$ It
is easy to see that $\mathcal Q(M)$ is a principal frame bundle of
$M$ with structure group $GL(n,H)Sp(1),$ it is also called a
$GL(n,H)Sp(1)$-structure on $M.$

Let $\pi :\mathcal Q(M)\longrightarrow M$ be the natural
projection. For each $u\in \mathcal Q(M)$ we consider the linear
isomorphisms $j(u)$ on $T_{\pi (u)}M$ defined by
$j(u)=uJ_{3}^{0}u^{-1}$. It is easy to see that $(j(u))^{2}=-id$.
For each point $p\in M$ we define $Z_{p}(M)=\{j(u): u\in \mathcal
Q(M), \pi (u)=p\}$. In other words, $Z_{p}(M)$ is the space of all
complex structures  in the tangent space $T_{p}M$ which are
compatible with the almost quaternionic structure on $M$.

We define the twistor space $Z$ of $(M,\mathcal Q)$ by setting
$Z=\bigcup_{p\in M}Z_{p}(M)$. Let $H_3$ be the stabilizer of
$J^{0}_3$ in the group $GL(n,H)Sp(1)$.  There is a bijective
correspondence between the symmetric space $GL(n,H)Sp(1)/H_3\cong
S^2= \{ (x,y,z)\in {\bf R}^3\ |\ x^2+y^2+z^2=1\}$ and $Z_{p}(M)$
for each $p \in M$. So we can consider $Z$ as the associated fibre
bundle of $\mathcal Q(M)$ with standard fibre $GL(n,H)Sp(1,{ \bf
R})/H_3.$ Hence, $\mathcal Q(M)$ is a principal fibre bundle over
$Z$ with structure group $H_3$ and projection $j$. We consider the
symmetric spaces $GL(n,H)Sp(1)/H_3.$ We have the following Cartan
decomposition $gl(n,H)\oplus sp(1)=h_3\oplus m_3$ where
$$h_3=\{A\in gl(n,H)\oplus sp(1): AJ^{0}_{3}=J^0_{3}A\}$$ is the
Lie algebra of $H_3$ and $m_3=\{A\in gl(n,H)\oplus sp(1):
AJ^0_{3}=-J^0_{3}A\}.$ It is clear that $m_3$ is generated by
$J^0_{1},$ $J^0_{2},$ i.e. $m_3=span\{J^0_{1},J^0_{2} \}.$  Hence,
if $A\in m_3$ then $J^0_{3}A\in m_3.$

Let $\nabla$ be a quaternionic connection on $M$, i.e. $\nabla$ is
a linear connection in the principal bundle $\mathcal Q(M)$ (see
e.g. \cite{15}). Note that we make no assumptions on the torsion
or on the curvature of $\nabla.$ Keeping in mind the similarity
with the 4-dimensional and  quaternionic geometry
\cite{2,11,S1,S3}, we use $\nabla$ to define two almost complex
structures $I^{\nabla}_1$ and $I^{\nabla}_2$ on the twistor space
$Z$. The construction of these almost complex structures depends
on the choice of the quaternionic connection $\nabla.$

We denote by $A^{\ast }$ (resp.  $B(\xi )$) the fundamental vector
field (resp. the standard horizontal vector field) on $\mathcal
Q(M)$ corresponding to $A\in gl(n,H)\oplus sp(1)$ (resp. $\xi \in
{\bf R^{4n}}$).

Let $u\in \mathcal Q(M)$ and $Q_{u}$ be the horizontal subspace of
the tangent space $T_{u}\mathcal Q(M)$ induced by $\nabla$ (see
e.g.  \cite{15} ). The vertical space i.e. the vector space
tangent to a fibre, is isomorhic to $${(gl(n,H)\oplus
sp(1))}^*_u=(h_3)^*_u\oplus (m_3)^*_u,$$ where $(h_3)^{\ast
}_{u}=\{A^{\ast }_{u}: A\in h_3\}, (m_3)^{\ast }_{u}=\{A^{\ast
}_{u}: A\in m_3\}.$

Hence, $T_{u}\mathcal Q(M)=(h_3)^{\ast }_{u}\oplus (m_3)^{\ast
}_{u}\oplus Q_{u}$.

For each $u\in \mathcal Q(M)$, we put $$V_{j(u)}=j_{\ast u}(
(m_3)^{\ast }_{u}), \quad H_{j(u)}=j_{\ast u}Q_{u}.$$ Thus we
obtain  vertical and horizontal distributions $V$ and $H$ on $Z$.
Since $\mathcal Q(M)$ is a principal fibre bundle over $Z$ with
structure group $H_3$  we have $Ker j_{\ast u}=(h_3)^{\ast }_{u}$.

Hence $V_{j(u)}=j_{\ast u}(m_3)^{\ast }_{u}$ and $j_{\ast
u|(m_3)^{\ast }_{u}\oplus Q_{u}}:(m_3)^{\ast }_{u}\oplus
Q_{u}\longrightarrow T_{j(u)}Z$ is an isomorphism.

We define two almost complex structures $I^{\nabla}_{1}$ and
$I^{\nabla}_{2}$ on $Z$ by
\begin{eqnarray}\label{2.2}
& &I^{\nabla}_{1}(j_{\ast u}A^{\ast })=j_{\ast u}(J^0_{3}A)^{\ast
}, \qquad I^{\nabla}_{2}(j_{\ast u}A^{\ast })=-j_{\ast
u}(J^0_{3}A)^{\ast }\\ & &I^{\nabla}_i(j_{\ast u}B(\xi ))=j_{\ast
u}B(J^0_{3}\xi ), \qquad i=1,2,\nonumber
\end{eqnarray}
for $A \in m_3, \xi \in {\bf R^{4n}}.$

For twistor bundles of 4-dimensional Riemannian manifolds the
almost complex structure $I^{\LC}_{1}$ is introduced in \cite{2}
and the almost complex structure $I^{\LC}_{2}$ is introduced in
\cite{11} in terms of the horizontal spaces of the Levi-Civita
connection. It is well known  that $I^{\LC}_{1}$ is integrable
exactly when the 4-manifold is anti-self-dual \cite{2}, while
$I^{\LC}_{2}$ is never integrable \cite{11}. The almost complex
structure $I^{\LC}_{1}$ and its integrability for QK spaces is
established in \cite{S1} and then generalized for quaternionic
manifold in \cite{S3}. The integrability of $I^{\nabla}_{1}$ in
the case of HKT and QKT manifolds is established in
\cite{HP3,HOP}, respectively. The non-integrability of
$I^{\LC}_{2}$ in the case of QK  is confirmed in \cite{DM} and the
non-integrability of $I^{\nabla}_1$ for QKT spaces  is done in
\cite{IM}.

\subsection{Dependence on the  quaternionic connection}
In this section we investigate when different  almost quaternionic
connections induce the same almost complex structure on the
twistor space over an almost quaternionic manofold.

Let $\nabla$ and $\nabla^{'}$ be two different  quaternionic
connections on an  almost quaternionic manifold $(M,\mathcal Q)$.
Then we have $$ \nabla^{'}_X=\nabla_X+S_X,\qquad X\in
\Gamma(TM),$$ where $S_X$ is a (1,1) tensor on $M$ and
$u^{-1}(S_X)u$ belongs to $gl(n, H)\oplus sp(1)$ for any $u\in
\mathcal Q(M)$. Thus we have the splitting
\begin{equation}\label{eq0}
S_X(Y)=S^0_X(Y)+s^1(X)J_1Y+s^2(X)J_2Y+s^3(X)J_3Y,
\end{equation}
where $X,Y\in \Gamma(TM)$, $s^i$ are 1-forms and $[S^0_X,J_i]=0,$
$i=1,2,3$.
\begin{pro}\label{t2.5}
Let $\nabla$ and $\nabla^{'}$ be two different
quaternionic connections on an  almost quaternionic manifold
$(M,\mathcal Q)$. The following conditions are equivalent:
\begin{enumerate}
\item[i).] The two almost complex structures $I_1^{\nabla}$ and
$I_1^{\nabla^{'}}$  on the twistor space $Z$ coincide.
\item[ii).] The 1-forms $s^1,s^2,s^3$ are related as follows
\begin{eqnarray*}
s^1(J_1X)=s^2(J_2X)=s^3(J_3X),\qquad X\in \Gamma(TM).
\end{eqnarray*}
\end{enumerate}
\end{pro}
\begin{proof}
We fix a point $J$ of the twistor space $Z$. We have
$J=a_1J_1+a_2J_2+a_3J_3$ with $a_1^2+a_2^2+a_3^2=1$.
Let $\pi: Z\longrightarrow M$ be the natural projection and
$x=\pi(J)$. The connection $\nabla$ induces a splitting of the
tangent space of $Z$ into vertical and horizontal components:
$T_JZ= V_J \oplus H_J$. Let $v$ and $h$ be the vertical and
horizontal projections corresponding to this splitting. Let $T_JZ=
{V^{'}}_J \oplus {H^{'}}_J$ be the splitting induced by
$\nabla^{'}$  with the projections $v^{'}$ and $h^{'}$,
respectively. It is easy to observe the following identities
\begin{eqnarray}\label{proj_ids} &
& v+h=1\nonumber \\ & &v^{'}+h^{'}=1 \\ & &vv^{'}=v^{'}
\nonumber\\ & & v^{'}+vh^{'}=v \nonumber
\end{eqnarray}
In fact, $V_J={V^{'}}_J$ and we may regard this space as a
subspace of  ${\mathcal Q}_x$. We have that
\begin{eqnarray*}
& & V_J=\{W\in {\mathcal Q}_x\ |\
WJ+JW=0\}=\{w_1J_1+w_2J_2+w_3J_3\ |\ w_1a_1+w_2a_2+w_3a_3=0\},
\end{eqnarray*}
where $J=a_1J_1+a_2J_2+a_3J_3$. It follows that for any
$W\in V_J$, $I_1^{\nabla}(W)=I_1^{\nabla^{'}}(W)=JW$. In general,
for any $W\in T_JZ$, we have
\begin{eqnarray}
& & I_1^{\nabla}(W)=J(vW)+(J\pi(W))^h \\ & &\nonumber
I_1^{\nabla^{'}}(W)=J(v^{'}W)+(J\pi(W))^{h^{'}},
\end{eqnarray}
where
$(.)^h$ (resp. $(.)^{h^{'}}$) denotes the horizontal lift on $Z$
of the corresponding vector field on $M$ with respect to $\nabla$
(resp. $\nabla^{'}$). Using (\ref{proj_ids}), we calculate that
\begin{gather}\label{coin}
v(I_1^{\nabla^{'}}W)=J(v^{'}W)+v(J\pi(W))^{h^{'}}=J((v-vh^{'})W)+v(J\pi(W))^{h^{'}}=\\\nonumber
v(I_1^{\nabla}W)-J(vh^{'}W)+v(J\pi(W))^{h^{'}}. \end{gather} We
investigate the equality
\begin{eqnarray}\label{eq1}
J(vh^{'}W)=v(J\pi(W))^{h^{'}},\qquad W\in T_JZ.
\end{eqnarray}
Take  $W=Y^{h^{'}}, Y\in \Gamma(TM)$ in (\ref{eq1}) to get
\begin{equation}\label{eq2}
J(vY^{h^{'}})=v(JY)^{h^{'}},\qquad Y\in T_xM
\end{equation}
Hence, (\ref{eq2}) is equivalent to
$I_1^{\nabla}=I_1^{\nabla^{'}}$ because of \eqref{coin}.

Let $(U,x_1,\dots,x_{4n})$ be a local coordinate system on $M$ and
let $Y=\sum Y^i\frac{\partial}{\partial x^i}$. The horizontal lift
of $Y$ with respect to $\nabla^{'}$ at the point $J\in Z$ is given
by
\begin{eqnarray}
Y^{h^{'}}_J=\sum_{i=1}^{4n}(Y^i\circ\pi)\frac{\partial}{\partial
x^i}-\sum_{s=1}^{3}a_s{\nabla^{'}}_YJ_s
\end{eqnarray}
We calculate
\begin{eqnarray}\label{eq3}
&
&v(JY)^{h^{'}}=(JY)^{h^{'}}-h(JY)^{h^{'}}=(JY)^{h^{'}}-(JY)^{h}=\\
& &\nonumber
=\sum_{s=1}^{3}a_s(-\nabla^{'}_{JY}J_s+\nabla_{JY}J_s)=-[S_{JY},J]
\end{eqnarray}
On the other hand, we have
\begin{eqnarray}\label{eq4}
J(vY^{h^{'}})=J(Y^{h^{'}}-Y^{h})=J\sum_{s=1}^{3}a_s(-\nabla^{'}_YJ_s+\nabla_YJ_s)=-J[S_Y,J]
\end{eqnarray}

Substitute (\ref{eq3}) and (\ref{eq4}) into (\ref{eq2}) to get
that $I_1^{\nabla}=I_1^{\nabla^{'}}$ is equivalent to the
condition
\begin{eqnarray}\label{eq5}
J[S_Y,J]=[S_{JY},J], \qquad Y\in \Gamma(TM), J\in Z.
\end{eqnarray}
Now, \eqref{eq5} easily leads to the equivalence of i) and ii).
\end{proof}
We note that \eqref{eq5} is discovered  in connection with the
coincidence of the almost complex structures generated by two
Oproiu connections in \cite{AMP1}.
\begin{co}\label{t2.72}
Let $\nabla$ and $\nabla^{'}$ be two different quaternionic
connections on an  almost quaternionic manifold $(M,\mathcal Q)$.
The following conditions are equivalent:
\begin{enumerate}
\item[i).] The two almost complex structures $I_2^{\nabla}$ and
$I_2^{\nabla^{'}}$  on the twistor space $Z$ coincide.
\item[ii).] The 1-forms $s^1,s^2,s^3$ vanish identically, $s_1=s_2=s_3=0$.
\end{enumerate}
\end{co}
\begin{proof}
It is sufficient to observe from the proof of
Proposition~\ref{t2.5} that $I_2^{\nabla}=I_2^{\nabla^{'}}$ is
equivalent to $J[S_Y,J]=-[S_{JY},J], \quad Y\in \Gamma(TM), J\in
Z$. The latter condition implies $s_1=s_2=s_3=0$.
\end{proof}
\begin{co}\label{t2.7}
 Let $\nabla$ and $\nabla^{'}$ be two different
quaternionic connections with torsion tensors $T^{\nabla^{'}}$ and
$T^{\nabla}$, respectively, on an  almost quaternionic manifold
$(M,\mathcal Q)$. The following conditions are equivalent:
\begin{enumerate}
\item[i).] The two almost complex structures $I_1^{\nabla}$ and
$I_1^{\nabla^{'}}$  on the twistor space $Z$ coincide.
\item[ii).] The $(0,2)_J$ part with respect to all $J\in\mathcal Q$ of the torsion $T^{\nabla}$ and
$T^{\nabla^{'}}$ coincides,
$$(T^{\nabla})^{0,2}_J=(T^{\nabla^{'}})^{0,2}_J.$$
\end{enumerate}
\end{co}
\begin{proof} Let $S=\nabla^{'}-\nabla$. Then we have
\begin{equation}\label{tr1}
T^{\nabla^{'}}(X,Y)=T^{\nabla}(X,Y)+S_X(Y)-S_Y(X).
\end{equation}
The $(0,2)_J$-part with respect to $J$ of \eqref{tr1} gives
\begin{gather}\label{eq7}
(T^{\nabla^{'}})^{0,2}_J-(T^{\nabla})^{0,2}_J=[S_{JX},J]Y -
J[S_X,J]Y-[S_{JY},J]X+ J[S_Y,J]X.
\end{gather}

For example, put  $J=J_3$ in (\ref{eq7}) and use the splitting
(\ref{eq0}) to obtain
\begin{gather}\label{eq8}
(T^{\nabla^{'}})^{0,2}_{J_3}-(T^{\nabla})^{0,2}_{J_3}=(s_1(J_3X)-s_2(X))J_2Y-(s_1(X)+s_2(J_3X))J_1Y\\\nonumber
-(s_1(J_3Y)-s_2(Y))J_2X+(s_1(Y)+s_2(J_3Y))J_1X.
\end{gather}
Hence we get $s_1(J_1X)=s_2(J_2X)$ is equivalent to
$(T^{\nabla^{'}})^{0,2}_{J_3}=(T^{\nabla})^{0,2}_{J_3}$.

Similarly, put $J=J_1$ in \eqref{eq7} and using \eqref{eq0} one
gets $s_3(J_3X)=s_2(J_2X)$ is equivalent to
$(T^{\nabla^{'}})^{0,2}_{J_1}=(T^{\nabla})^{0,2}_{J_1}$.  Hence,
$s_1(J_1X)=s_2(J_2X)=s_3(J_3X)$ is equivalent to
$(T^{\nabla^{'}})^{0,2}_{J_3}=(T^{\nabla})^{0,2}_{J_3}$ and
$(T^{\nabla^{'}})^{0,2}_{J_1}=(T^{\nabla})^{0,2}_{J_1}$. The
latter conditions imply
$(T^{\nabla^{'}})^{0,2}_{J_2}=(T^{\nabla})^{0,2}_{J_2}$  since the
formula (3.4.4) in \cite{AM}) expressing $N_{J_3}$ by $N_{J_1}$
and $N_{J_2}$ holds for the $(0,2)_{J_a}$-part $T^{\nabla}_{J_a},
a =1,2,3$, of the torsion. It is easy to see that the general
formula (6) in \cite{AMP1} which expresses $N_J$ in terms of
$N_{J_1}, N_{J_2}, N_{J_3}$ holds also for the $(0,2)_J$-part of
any tensor from $\Lambda^2(TM)$. Applying this formula to the
$(0,2)_J$-part of the torsion, we conclude that
$(T^{\nabla^{'}})^{0,2}_J=(T^{\nabla})^{0,2}_J$ is equivalent to
$s_1(J_1X)=s_2(J_2X)=s_3(J_3X)$. Proposition~\ref{t2.5} completes
now the proof.
\end{proof}

\subsection{Integrability}

In this section we investigate conditions on the almost
quaternionic connection $\nabla$ which imply the integrability of
the almost complex structure $I^{\nabla}_1$ on $Z$. We also  show
that $I^{\nabla}_2$ is never  integrable i.e. for any choice of
the almost quaternionic connection  $\nabla$ it has non-vanishing
Nijenhuis tensor.

We denote by $IN_i,  \quad i=1,2$ the Nijenhuis tensors of $I_i$
and  recall that
\begin{gather*}
IN_i(U,W)=[I_iU,I_iW]-[U,W]-I_i[I_iU,W]-I_i[U,I_iW], \quad U,W\in
\Gamma(TZ).
\end{gather*}
\begin{pro}\label{lint}
Let $\nabla$ be a quaternionic connection on an
almost quaternionic manifold $(M,\mathcal Q)$ with torsion tensor
$T^{\nabla}$.  The following conditions are equivalent:
\begin{enumerate}
\item[i).] The almost complex structure $I^{\nabla}_1$ on the twistor
space $Z$ of $(M,\mathcal Q)$ is integrable.
\item[ii).] The $(0,2)_J$-part $(T^{\nabla})^{0,2}_J$ of the torsion with respect
to all $J\in \mathcal Q$ vanishes,  and  the (2,0)+(0,2) parts of
the Ricci 2-forms with respect to an admissible basis
$J_1,J_2,J_3$ of $\mathcal Q$ coincide, i.e. the next conditions
hold
\begin{gather}\label{ltor}
 (T^{\nabla})^{0,2}_J=0, J\in \mathcal Q,\\
\label{idric}  \rho_a(J_cX,J_cY) - \rho_a(X,Y)+\rho_b(J_cX,Y) +
\rho_b(X,J_cY)=0.
\end{gather}
\end{enumerate}
\end{pro}
\begin{proof} Let $J_1,J_2,J_3$ be an admissible basis
 of the almost quaternionic structure $\mathcal
Q$.

Let $hor$ be the natural projection $T_uQ\longrightarrow
(m_3)^*_u\oplus Q_u,$ with $ker(hor)=(h_3)^*_u$. We define a
tensor field $I^{'}_1$ on $\mathcal Q(M)$ by
\begin{eqnarray*}\nonumber
& &I^{'}_{1}(U) \in  (m_3)^*_u\oplus Q_u, \\ &
&j_{*u}(I^{'}_1(U))=I_1(j_{*u}U), \qquad U\in T_uP.\nonumber
\end{eqnarray*}
For any $U, W\in \Gamma(T\mathcal Q(M))$ we define
$$IN^{'}_1(U,W)=hor[I^{'}_1U,I^{'}_1W]-hor[hor U, hor W]-
I^{'}_1[I^{'}_1U,hor W]-I^{'}_1[hor U,I^{'}_1W]$$

It is easy to check that $IN^{'}_1$ is a  a tensor field on
$\mathcal Q(M)$. We also observe that
\begin{equation}\label{in1}
j_{*u}(IN^{'}_1(U,W))=IN_1(j_{*u}U,j_{*u}W),
\qquad U,W\in T_u\mathcal Q(M)
\end{equation}
Let $A,B\in m_3$ and $\xi, \eta\in {\bf R}^{4n}$. Using the well
known general commutation relations  among  the fundamental vector
fields and standard horizontal vector fields on the principal
bundle $\mathcal Q(M)$ (see e.g. \cite{15}), we calculate taking
into account \eqref{in1} that
\begin{eqnarray}
& &IN_1(j_{*u}(A^*_u),j_{*u}(B^*_u))=0. \nonumber\\ &
&IN_1(j_{*u}(A^*_u),j_{*u}(B(\xi)_u))= 0.\nonumber\\\label{n1} &
&[IN_1(j_{*u}(B(\xi)_u)),j_{*u}(B(\eta)_u))]_{H} =
\\ & &\qquad j_{*u}
(B(-\Theta(B(J_3^0\xi),B(J_3^0\eta))+\Theta(B(\xi),B(\eta))\nonumber\\
& &\qquad
+J_3^0\Theta(B(J_3^0\xi),B(\eta))+J_3^0\Theta(B(\xi),B(J_3^0\eta)))_u).
\nonumber\\\label{n2}
& &[IN_1(j_{*u}(B(\xi)_u)),j_{*u}(B(\eta)_u))]_{V}=\\& &\qquad
\{-\rho_1(B(J_3^0\xi),B(J_3^0\eta)) + \rho_1(B(\xi),B(\eta)
)\nonumber\\ & &\qquad -\rho_2(B(J_3^0\xi),B(\eta)) -
\rho_2(B(\xi),B(J_3^0\eta))\} j_{*u}(J_1^0)\nonumber\\ &
&\nonumber\qquad +\{-\rho_2(B(J_3^0\xi),B(J_3^0\eta)) +
\rho_2(B(\xi),B(\eta))\\ & &\qquad +\rho_1(B(J_3^0\xi),B(\eta)) +
\rho_1(B(\xi),B(J_3^0\eta)))\}
j_{*u}(J_2^0).\nonumber\\\label{n33}& &
IN_2(j_{*u}(A^*_u),j_{*u}(B(\xi)_u))= -4j_{*u}(B(A\xi)_u)\not=0.
\end{eqnarray}

Take $X=u(\xi), Y=u(\eta)$, we see that \eqref{n1} and \eqref{n2}
are equivalent to
\begin{gather}\label{tor}
(T^{\nabla})^{0,2}_{J_3}=T^{\nabla}(J_3X,J_3Y)-T^{\nabla}(X,Y)-J_3T^{\nabla}(J_3X,Y)-J_3T^{\nabla}(X,J_3Y)=0,\\
\label{curv}
(\Pi^{\nabla})^{0,2}_{J_3}=\rho^{\nabla}_1(J_3X,J_3Y)-\rho^{\nabla}_1(X,Y)+\rho^{\nabla}_2(J_3X,Y)+
\rho^{\nabla}_2(X,J_3Y)=0,
\end{gather}
respectively.

Similarly, we get
$(\Pi^{\nabla})^{0,2}_{J_3}=(\Pi^{\nabla})^{0,2}_{J_2}=(\Pi^{\nabla})^{0,2}_{J_1}=0$
and
$(T^{\nabla})^{0,2}_{J_3}=(T^{\nabla})^{0,2}_{J_2}=(T^{\nabla})^{0,2}_{J_1}=0$.
The first equalities imply that \eqref{n2} is equivalent to
\eqref{idric}. We apply to the second equalities the same
arguments as in the proof of Corollary~\ref{t2.7}, i.e. formula
(6) in \cite{AM}, to derive that \eqref{n1} is equivalent to
$(T^{\nabla})^{0,2}_{J}=0$ for all local $J\in\mathcal Q$.
\end{proof}
The equation \eqref{n33}  in the proof of Proposition~\ref{lint}
yield
\begin{co}
Let $\nabla$ be an almost quaternionic connection on an almost
quaternionic manifold $(M,\mathcal Q)$ with torsion tensor
$T^{\nabla}$. Then the almost complex structure $I^{\nabla}_2$ on
the twistor space $Z$ of $(M,\mathcal Q)$ is never integrable.
\end{co}
In the 4-dimensional case we derive
\begin{thm}\label{four}
Let $(M^4,g)$ be a 4-dimensional Riemannian
manifold with a Riemasnnian metric $g$ and let $\mathcal Q$ be the
quaternionic structure corresponding to the conformal class
generated by $g$ with a local basis $J_1,J_2,J_3$. Then the
following conditions are equivalent
\begin{enumerate}
\item[i).] The  metric $g$ is anti-self-dual.
\item[ii).] The Ricci forms $\rho_a^g$ of the Levi-Civita
connection $\LC$ satisfy \eqref{idric}, i.e.
\begin{equation*} \rho^g_a(J_cX,J_cY) -\rho^g_a(X,Y)+\rho^g_b(J_cX,Y)
+ \rho^g_b(X,J_cY)=0.
\end{equation*}
\item[iii).]  The torsion condition
\eqref{ltor}  for a linear connection $\nabla$ always implies the
curvature condition \eqref{idric}.
\end{enumerate}
\end{thm}
\begin{proof}
The proof is a direct consequence of Proposition~\ref{lint},
Corolarry~\ref{t2.7} and the result in \cite{2} which states that
the almost complex structure $I_1^{\LC}$ is integrable exactly
when the conformal structure generetaed by $g$ is anti-self-dual.
\end{proof}
In higher dimensions, the curvature condition \eqref{idric} is a
consequence of the torsion condition \eqref{ltor} in the sense of
the next
\begin{thm}\label{t2.1}
Let $\nabla$ be a quaternionic connection on an
almost quaternionic 4n-dimensional $n\ge 2$ manifold $(M,\mathcal
Q)$ with torsion tensor $T^{\nabla}$. Then the following
conditions are equivalent:
\begin{enumerate}
\item[i).] The almost complex structure $I^{\nabla}_1$ on the twistor
space $Z$ of $(M,\mathcal Q)$ is integrable.
\item[ii).] The $(0,2)_J$-part $(T^{\nabla})^{0,2}_J$ of the torsion with respect
to all $J\in \mathcal Q$ vanishes,
\\\centerline{$(T^{\nabla})^{0,2}_J=0, J\in \mathcal Q$.}
\end{enumerate}
\end{thm}
\begin{proof} Suppose i) holds. Then ii) follows from Proposition~\ref{lint}.

For the converse, \eqref{ltor} and the fact that  the connection
$\nabla$ is a quaternionic connection, $\nabla\in\Delta(\mathcal
Q)$,  yield the next expression for the Nijenhuis tensor $N_J$ of
any local $J\in\mathcal Q$,
\begin{equation}\label{pq0}
N_J(X,Y)\in
span\{J_1X,J_1Y,J_2X,J_2Y,J_3X,J_3Y\},
\end{equation}
where $J_1,J_2,J_3$ is an admissible local basis of $\mathcal Q$.

To prove that ii) implies the integrability of $I_1^{\nabla}$, we
apply the result in \cite{AM} (see also \cite{AMP1}) which states
that an almost quaternionic 4n-maniofold $(n\ge 2)$ is
quaternionic if and only if the three Nijenhuis tensors
$N_1,N_2,N_3$ satisfy the condition
\begin{equation}\label{pq}
(N_1(X,Y)+N_2(X,Y)+N_3(X,Y))\in
span\{J_1X,J_1Y,J_2X,J_2Y,J_3X,J_3Y\}.
\end{equation}
Clearly, \eqref{pq} follows from \eqref{pq0} which shows that the
almost quaternionic 4n-manifold ($n\ge2$) $(M,\mathcal Q)$ is a
quaternionic manifold. Let  $\nabla^0$ be a torsion-free
quaternionic  connection on $(M,\mathcal Q)$. Then the almost
complex structure $I_1^{\nabla^0}$ on the twistor space $Z$ is
integrable \cite{S3} and $I_1^{\nabla}=I_1^{\nabla^0}$ due to
Corolarry~\ref{t2.7}.

Hence, the equivalence between i) and ii) is established, which
completes the proof.
\end{proof}
From the proof of Proposition~\ref{lint} and Theorem~\ref{t2.1},
we easily derive
\begin{co}\label{cur}
Let $\nabla$ be an almost quaternionic connection on an
$4n$-dimensional ($n\ge 2$) almost quaternionic manifold
$(M,\mathcal Q)$ with torsion tensor $T^{\nabla}$. Then the
torsion condition \eqref{ltor} implies the curvature condition
\eqref{idric}.
\end{co}

We note that Corollary~\ref{cur} generalizes the same statement
derived in the case of QKT-connection on QKT manifolds in
\cite{IM}.

Theorem~\ref{t2.1} and Corollary~\ref{t2.7}  imply
\begin{co}\label{t2.8}
Let $(M,\mathcal Q)$ be an almost quaternionic
manifold. Among the all almost complex structures $I_1^{\nabla},
\nabla\in\Delta(\mathcal Q)$ on the twistor space $Z$ at most one
is integrable.
\end{co}
The proof  of the next theorem follows directly from the proof of
Theorem~\ref{t2.1}, Theorem~\ref{four} and Corollary~\ref{t2.8}.
\begin{thm}\label{main}
Let $(M,\mathcal Q)$ be an almost quaternionic $4n$-manifold. The
next two conditions are equivalent:
\begin{enumerate}
\item[1).]  Either $(M,\mathcal Q)$ is a quaternionic manifold (if
$n\ge 2$) or $(M,\mathcal Q=[g])$ is ant-self dual for $n=1$.
\item[2).] There exists an integrable almost complex structure
$I^{\nabla}_1$ on the twistor space $Z$ which does not depend on
the quaternionic connection $\nabla$.
\end{enumerate}
\end{thm}

\section{Quaternionic K\"ahler manifolds with torsion }
An almost quaternionic Hermitian manifold $(M,\mathcal Q,g)$ is
called quaternionic K\"ahler with torsion (QKT) if there exists a
an almost quaternionic Hermitian connection
$\nabla^T\in\Delta({\mathcal Q})$ whose torsion tensor $T$ is a
3-form which is (1,2)+(2,1) with respect to each $J_a$ \cite{HOP},
i.e. the tensor $T(X,Y,Z):=g(T(X,Y),Z)$ is totally skew-symmetric
and satisfies the conditions
\begin{gather*}
 T(X,Y,Z)=T(J_aX,J_aY,Z) + T(J_aX,Y,J_aZ)
+T(X,J_aY,J_aZ), \quad a=1,2,3.
\end{gather*}

We recall that each QKT is a quaternionic manifold due to an
observation made in \cite{I1}. The condition on the torsion
implies that the (0,2)-part of the torsion of a QKT connection
vanishes. Applying Theorem~\ref{t2.1}, we obtain
\begin{thm}\label{qkt}
Let $(M,\mathcal Q,\nabla^T)$ be a QKT and
$\nabla^0\in\Delta({\mathcal Q})$ be a torsion-free quaternionic
connection. Then  the complex structure $I_1^{\nabla^T}$  on the
twistor space $Z$ constructed in \cite{HOP}    coincides with the
complex structure $I_1^{\nabla^0}$ constructed in \cite{S3}.
\end{thm}

\bibliographystyle{hamsplain}

\providecommand{\bysame}{\leavevmode\hbox
to3em{\hrulefill}\thinspace}

\end{document}